\documentclass[12pt]{article}

\usepackage{pstricks,pst-node,graphicx}
\usepackage{amssymb,amsmath}

\newcommand{\be}{\begin{equation}}
\newcommand{\ee}{\end{equation}}
\newcommand{\eea}{\end{eqnarray}}
\newcommand{\bea}{\begin{eqnarray}}
\newcommand{\xin}{\xi_N}
\newcommand{\szn}{Z_{\xin}}
\newcommand{\tzn}{\tilde Z_{N}}
\newcommand{\zn}{Z_{N}}
\newcommand{\numb}{n_b+n_r}
\newcommand{\leng}{n_b+n_r+n_{br}}
\newcommand{\tnumb}{n_b+n_r+\gamma_b+\gamma_r}

\newcommand{\uvw}{(u,v,w)}
\newcommand{\uno}{\vert_{(u,v,w)=(1,1,1)}}
\newcommand  \nt {{\noindent}}

\newcommand{\rv}{{r.v.}}
\newcounter{secnum}[section]
\begin{document}


\title{{{\bf Central Limit Theorem for Coloured Hard-Dimers}}
\footnotetext{\newline{\bf Key words and phrases}: Coloured hard-dimers, generating function, probability distribution, central limit theorem
\newline
{\it Mathematics Subject Classification}: 60F05, 05A15, 60C05}}

\author{{\bf Maria Simonetta Bernabei} \\ and \\
{\bf Horst Thaler} \\[1ex] Department of Mathematics and Informatics,  \\
University of Camerino, \\
Via Madonna delle Carceri 9,
I--62032, Camerino (MC), Italy;\\
{\small simona.bernabei@unicam.it, horst.thaler@unicam.it}}

\date{}

\maketitle

\abstract{{Using an averaged generating function for coloured
hard-dimers, some random variables of interest are studied. The main result lies in the fact that all their probability distributions obey a central limit theorem.

}}

\section{Introduction}
 \setcounter{secnum}{\value{section}
 \setcounter{equation}{0}
 \renewcommand{\theequation}{\mbox{\arabic{secnum}.\arabic{equation}}}}
\newtheorem{remark}{Remark}[section]
\newtheorem{lemma}{Lemma}[section]
\newtheorem{proposition}{Proposition}[section]
\newtheorem{theorem}{Theorem}[section]
\newtheorem{cor}{Corollary}[section]

In the literature coloured hard-dimers are applied in the framework of causally triangulated $(2+1)$-dimensional quantum gravity. It was proved in \cite{BeLoZa}, by using special triangulations of spacetime, that the generating function of the one step propagator depends on that one of coloured hard-dimers. For any configuration $\xin$, of length $N$, of blue and red sites, a coloured hard-dimer is a sequence of blue and red dimers, that satisfy the ``hardness" condition, i.e. they can not intersect. A dimer is an edge connecting two nearest sites of the same colour.

In the present paper we consider the ``averaged" generating function for coloured hard-dimers, that is, the mean of generating functions over all the configurations $\xin$, with $N$ fixed. In
\cite{BeTh} we have found an explicit formula for it together with estimates from above and from below, that are  both exponential, for large $N$. In the following we study the probability distribution associated with the averaged generating function. Then we analyze the probability distributions corresponding to some random variable (\rv) of interest. In particular the number of dimers and the total length of them. It turns out that the role played by the r.v. that measures the total number of dimers and single points (i.e. sites not occupied by dimers) is very important. We prove that the r.v. total length of dimers is binomial with parameters $N-1$ and $\frac 13$. Moreover we see that, even though the number of dimers has an unknown probability distribution, we are able to estimate its mean and variance asymptotically, by using some recursive formulas that relate its moments with that ones resulting from the dimers' length.
Although the dimers' number distribution is not binomial, its variance is of order $N$, as $N \to \infty$, as in the binomial case.

The main result of the present article is a local Central Limit Theorem (C.L.T.), for $N$ large enough, for the joint probability distribution corresponding to the number of dimers and that one of dimers and single points. The limit distribution is a bivariate gaussian distribution with correlation coefficient equal to $-{1 \over \sqrt{3}}$. Hence a C.L.T. holds, as $N \to \infty$, also for the marginal probability distribution related to the dimers' number.

The paper is organized as follows. In section 2 we define the probability distribution associated with coloured hard-dimers, through the averaged generating function, and find an exact expression for its normalizing constant $C_N=\left ( {3 \over 2}\right)^{N-1}$. Moreover we recognize the right probability distribution for the length of dimers. In section 3 we calculate the first two moments of the dimers' number. Finally in section 4 we prove a C.L.T. for the dimers' number.

\section{ Coloured hard-dimers and probability distributions}
 \setcounter{secnum}{\value{section}
 \setcounter{equation}{0}
 \renewcommand{\theequation}{\mbox{\arabic{secnum}.\arabic{equation}}}}

Given a sequence $\xi_N$ of length $N$ of blue and red sites on the one-dimensional lattice $\mathbb{Z}$, one defines a dimer to be an edge connecting two nearest sites of the same colour, that characterizes the dimer colour. A sequence of coloured and non overlapping dimers in turn yields a ``coloured hard-dimer". In Fig.1 an example of a coloured hard-dimer is given.

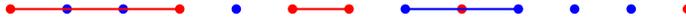
\begin{figure}[h]
\begin{center}
\psset{xunit=1.5cm,yunit=0.75cm}
\begin{pspicture}(0,0)(7,2)
\psdots[dotstyle=*,linecolor=red](0.5,0.5)(2.0,0.5)(3,0.5)(3.5,0.5)(4.5,0.5)(6.5,0.5)
\psdots[dotstyle=*,linecolor=blue](1,0.5)(1.5,0.5)(2.5,0.5)(4,0.5)(5,0.5)(5.5,0.5)(6,0.5)
\psline[linecolor=red](0.5,0.5)(2,0.5)
\psline[linecolor=red](3,0.5)(3.5,0.5)
\psline[linecolor=blue](4,0.5)(5,0.5)
\end{pspicture}
\caption{\label{hdimer} A hard-dimer, $N=13$}
\end{center}
\end{figure}

As described in \cite{BeLoZa} and \cite{BeTh} one introduces the generating function associated to coloured hard-dimers on $\xin$
$$
\szn(u,v,w) = \sum_{D} u^{n_b(D)} v^{n_r(D)} w^{n_{br}(D)},\quad u,v,w\in (0,\infty),
$$
where $D$ is a hard-dimer on $\xin$, $n_b(D)$ and $n_r(D)$ indicate the number
of blue and red dimers respectively on  $D$, and $n_{br}(D)$ the total number of sites within each dimer (this means the sites having a colour different from the colour of dimers containing them).
Moreover, let us define by $\gamma_b(D)$ and $\gamma_r(D)$ the number of blue and red sites of $D$ respectively, not occupied by dimers (``single points"). Then the following constraint
\be
\label{nr1}
2 n_b(D) + 2n_r(D) + n_{br}(D)+  \gamma_b(D) +\gamma_r(D)=N
\ee
holds.
\noindent
In the above example $n_b(D)=1,\,  n_r(D)=2,\,  n_{br}(D)=3,\,  \gamma_b(D)=3$ and $\gamma_r(D)=1$.

In \cite{BeTh} we studied the average of $\szn$ over all the configurations $\xin$, i.e.

\be
\label{nr2}
\tzn(u,v,w) = {1 \over {2^{N}}} \sum_{\xin} \szn (u,v,w)
\ee
and found estimates from above and from below for $\tzn$. In order to obtain them we proved an explicit formula for the mean $\tzn$, by using combinatorial tools

\be
\label{nr3}
\tzn \uvw =
1+  \sum\limits_{t=1}^N \sum\limits_{s=1}^{[{\frac t2}]} {N-t+s \choose s} {t-s-1 \choose s-1} \left ({{u+v}\over 4}\right )^s \left( {w \over 2} \right )^{t-2s}
 \ee
($[\cdot]$ denotes the integer part).

Now, let us define a family of probability spaces $(\Omega_N,{\cal F}_N,P_N)_{N\in \mathbb{N}_{\geq 1}}$. We choose $\Omega_N$ to be the set of all different hard-dimer configurations, where by hard-dimer configuration we mean a particular sequence $\xin$ of length $N$ together with a particular hard-dimer on $\xin$. The $\sigma$-algebra  ${\cal F}_N$ is the set of all subsets of $\Omega_N$ and for $P_N$ we take the probability measure which is distributed uniformly on $\Omega_N$. Normalizing the function $\tzn$ by $\tzn \uno$ then just gives the joint generating function of the random variables $n_b,n_r,n_{br}$, defined on $(\Omega_N,{\cal F}_N,P_N)$, which count, for each hard-dimer configuration, the number of blue, red dimers and the total number of sites within dimers, respectively.

The random variable $\numb$ (characterized by the index $s$ in formula (\ref{nr3})) gives the number of dimers for each hard-dimer configuration. Correspondingly, the random variable  $N-  \gamma_b- \gamma_r$ (indexed by $t$ in (\ref{nr3})) denotes the number of sites occupied by dimers. The above formula (\ref{nr3}) is obtained by fixing first the number of blue dimers ($n_b$), red ones ($n_r$) and single points ($\gamma_b$ and $\gamma_r$), i.e., $\tnumb$, indexed by $N-t+s$, without assigning any length to dimers, even though the total length of dimers, defined as $\leng$, (indexed by $t-s$)
is given because of (\ref{nr1}). This gives a factor
$$
{{(\tnumb)!}\over{n_b! \;  n_r! \; \gamma_b! \gamma_r!}}
$$
Then, in the non-trivial case where $n_b+n_r \neq 0$, one can assign a length to each dimer, taking into account that the total length of them is given by $\leng$, contributing another factor
${ {\leng - 1}\choose {\numb-1}}$. By summing over $\gamma_b, n_b$ and then over $s$ and $t$, we get the expression (\ref{nr3}). See \cite{BeTh} for the details.

In the present paper we want to analyze the random variables $\numb$ and $\leng$, more precisely their probability distribution. We prove that the random variable $\leng$ is binomial with parameters $N-1$ and ${1 \over 3}$ and, hence, a Central Limit Theorem (C.L.T.) holds, for large $N$ (De Moivre-Laplace's Theorem). In the case of $\numb$ we see that the probability distribution is unknown, but, fixing the random variable $\tnumb$, which in turn is binomial with parameters $N-1$ and ${2 \over 3}$, because of (\ref{nr1}), the conditional distribution of $\numb$ is hypergeometric. Moreover, for $N$ large enough, a C.L.T. holds for the joint probability distribution of $\numb$ and $\tnumb$ and, hence, also for the distribution of $\numb$. For the proof of these results the random variable $\tnumb$ plays a very important part.

Evaluating the averaged generating function $\tzn$ at the point $u=v=w=1$ we derive the normalizing constant of the probability measure $P_N$ associated to coloured hard-dimers
\be
\label{nr4}
C_N = \tzn \uno =
1+  \sum\limits_{t=1}^N \sum\limits_{s=1}^{[{\frac t2}]} {N-t+s \choose s} {t-s-1 \choose s-1} {{1}\over {2^{t-s}}}
\ee
This also shows that the joint probability distribution $\tilde{P}_N$ related to the r.v.s $n_b+n_r$ and $2(n_b+n_r)+n_{br}$, more precisely
\be
  \label{nr6}
  {\tilde P}_N(s,t)= P_N\left (\numb =s; 2 (\numb) +n_{br}=t \right )
  \ee
is given by
\be
\label{nr5}
 {\tilde P}_N(s,t) = \left \{ \begin{array} {ll}
{1 \over {C_N}}  & \textnormal{for $s=t=0$} \\
{1 \over {C_N}} {N-t+s \choose s} {t-s-1 \choose s-1} {{1}\over {2^{t-s}}}
 & \textnormal{otherwise}
 \end{array} \right.
 \ee
\nt
The main result of this section is an explicit formula for the normalizing constant $C_N$, that holds for any $N$. It is obtained by using combinatoric arguments.

  \begin{theorem}
  For any $N$, the following formula
  \be
  \label{nr7}
  C_N =\left ( {3 \over 2 } \right )^{N-1}
  \ee
   holds.
  \end{theorem}

\nt
{\bf Proof:}
Consider the following change of variables $t^{\prime} = N-t, k=N-t+s$, so that (\ref{nr4}) becomes (note that after changing the variables, we shall rename $t'$ again by $t$)

\be
\label{nr8}
C_N =
1+ \left ( \sum\limits_{k=1}^{[\frac N2]} \; \sum\limits_{t=0}^{k-1}+ \sum\limits_{k=[{\frac N2}] +1}^{N-1} \; \sum\limits_{t=2k-N}^{k-1}\right ){k \choose t} {N-k-1 \choose k-t-1} {{1}\over {2^{N-k}}}
\ee
Note that in formula (\ref{nr8}) only the combinatorial coefficients depend on $t$. Moreover, the combinatorial coefficients of the first sum in (\ref{nr8})  ${k \choose t} {N-k-1 \choose k-t-1}$ yield a non-normalized hypergeometric distribution with parameters $N-1$ (population size), $k$ (number of successes in the population) and $k-1$ (sample size). Therefore summing over $t$ we get
\be
\label{nr9}
\sum\limits_{k=1}^{[{\frac N2}]} {{1}\over {2^{N-k}}}  \sum\limits_{t=0}^{k-1} {k \choose t} {N-k-1 \choose k-t-1} = \sum\limits_{k=1}^{[\frac N2]} {{1}\over {2^{N-k}}} {N-1 \choose k-1}
\ee
  Analogously, performing the change of variable $t^{\prime}=t-2k+N$, we get in the second sum of (\ref{nr8})  again a non-normalized hypergeometric distribution with parameters $N-1, N-k-1, N-k-1$. Summing again over $t$ we get

$$ \sum\limits_{k=[{\frac N2}] +1}^{N-1}  {{1}\over {2^{N-k}}} \sum\limits_{t=2k-N}^{k-1} {k \choose t} {N-k-1 \choose k-t-1} =$$
 $$ \sum\limits_{k=[{\frac N2}] +1}^{N-1}  {{1}\over {2^{N-k}}}   \sum\limits_{t=0}^{N-k-1}  {k \choose N-k-t} {N-k-1 \choose t} = $$
\be
\label{nr10}
\sum\limits_{k=[{\frac N2}] +1}^{N-1}  {{1}\over {2^{N-k}}}  {N-1 \choose k-1}
\ee
with the convention that the binomial coefficient ${n \choose k}=0$ if $k>n$. Putting together the last terms in (\ref{nr9}) and (\ref{nr10}) we get a binomial formula
\be
\label{nr11}
\sum\limits_{k=1}^{N-1}  {{1}\over {2^{N-k}}}  {N-1 \choose k-1}
={{1}\over {2^{N-1}}}
\sum\limits_{k=0}^{N-1}   {N-1 \choose k} 2^k -1= \left ( {3 \over 2}\right )^{N-1}-1
\ee

\hfill{$\Box$}

An important Corollary of the previous Theorem is that the probability distribution of the r.v. $\tnumb$ (total number of dimers and single points) related to the index $k$ in (\ref{nr8}) has a binomial distribution with parameters $N-1$ and $\frac 23$.

\medskip
\nt
\begin{cor} The probability distribution of the r.v. $\tnumb$ is a binomial distribution with parameters $N-1$ and $\frac 23$.
\end{cor}

\nt
{\bf Proof: } Performing the same change of variables as for (\ref{nr8}) and summing the probability distribution with respect to $t$, as in (\ref{nr8}) we get
$$P_N(\tnumb=k) = \left( {2 \over 3}\right )^{N-1} {N-1 \choose k-1}  { 1 \over {2^{N-k}}}=
$$
$$
{N-1 \choose k-1}  \left( {2 \over 3}\right )^{k-1} \left( {1 \over 3}\right )^{N-k}$$
for $1 \leq k\leq N-1$. For $k=N$ it holds $\numb=0 \Rightarrow \gamma_b + \gamma_r =N$. Hence $$P_N(\tnumb =N) =P_N(\numb=0)= {1 \over {C_N}} =  \left( {2 \over 3}\right )^{N-1}$$ by (\ref{nr5}).
\hfill{$\Box$}
\medskip

\begin{cor} For large $N$ a C.L.T. (De Moivre-Laplace) for the r.v. $\tnumb$ holds
$$P_N(\tnumb=k)= {1 \over {\sqrt{ {\frac 49} \pi (N-1)}}}
e^{-{{\left (  k - {\frac 23} (N-1)\right )^2} \over {{\frac 49}(N-1)}}}
 \left (1 + r_k(N) \right )$$
where $\lim_{N \to \infty} r_k (N) =0$, uniformly with respect to $k$, for any $k$ such that
$y={{\left [  k - {\frac 23} (N-1) \right ]} / {\sqrt{{\frac 29}(N-1)}}}$ belongs to a finite interval $(-A, A)$.
\end{cor}

Analogously one can find a similar result for the r.v. $\leng$ (total length of dimers).

\begin{cor} The probability distribution of the r.v. $\leng$ is a  binomial distribution with parameters $N-1$ and $\frac 13$.
\end{cor}

\nt
{\bf Proof: } Taking into account (\ref{nr1}) one has
$$P_N(\leng=h) = P_N(\tnumb =N-h)=
$$
$$
{N-1 \choose h}  \left( {1 \over 3}\right )^{h} \left( {2 \over 3}\right )^{N-1-h}$$
for $0 \leq h\leq N-1$.
\hfill{$\Box$}

\begin{cor} For large $N$ a C.L.T. (De Moivre-Laplace) for the r.v. $\leng$ holds
$$P_N(\leng=h)= {1 \over {\sqrt{ {\frac 49} \pi (N-1)}}}
e^{-
{
{\left (  h - {\frac 13} (N-1)\right )^2} \over {{\frac 49}(N-1)}
}
} \left (1 + r'_h(N) \right )$$
where $\lim_{N \to \infty}r'_h(N) =0$, uniformly with respect to $h$, for any $h$ such that
$z={{\left [ h - {\frac 13} (N-1)\right ]}/ {\sqrt{{\frac 29}(N-1)}}}$ belongs to a finite interval $(-B, B)$.
\end{cor}

\section{Number of dimers: moments}
 \setcounter{secnum}{\value{section}
 \setcounter{equation}{0}
 \renewcommand{\theequation}{\mbox{\arabic{secnum}.\arabic{equation}}}}
In the present section we investigate the distribution of the random variable $\numb$, more precisely, we calculate the first two moments of it, by using recursive asymptotic formulas that depend on the first two moments of the binomial r.v. $\leng$, studied in the previous section. Starting from the averaged generating function $\tzn$, defined in (\ref{nr3}), we  rescale it by the normalizing constant $C_N$, calculated in  the previous section
$$\zn\uvw={{\tzn\uvw}\over {C_N}}= $$
\be
\label{nr31}
\left ( {2 \over 3}\right )^{N-1}  \sum\limits_{t=1}^N \sum\limits_{s=1}^{[{\frac t2}]} {N-t+s \choose s} {t-s-1 \choose s-1} {1 \over {2^{t-s}}}\left ({{u+v}\over 2}\right )^s w^{t-2s}
\ee
Therefore we calculate the moments of $\numb$ through its derivatives (see \cite{GrSt}). By symmetry of the variables $u$ and $v$ in (\ref{nr31}) we have
$$ E_N(n_b)= {\partial \over {\partial u}} \zn\uvw \uno = $$
$$
{\partial \over {\partial v}} \zn \uvw \uno = E_N(n_r)=$$
$$
\left ( {2 \over 3}\right )^{N-1}  \sum\limits_{t=1}^N \sum\limits_{s=1}^{[{\frac t2}]} {{s /2} \over {2^{t-s}}}  {N-t+s \choose s} {t-s-1 \choose s-1}
$$
We indicate by $E_N$ the mean of random variables with respect to the probability measure $P_N$.
Therefore
$$
E_N(\numb) =
\left ( {2 \over 3}\right )^{N-1}  \sum\limits_{t=1}^N \sum\limits_{s=1}^{[{\frac t2}]}{s \over {2^{t-s}}}  {N-t+s \choose s} {t-s-1 \choose s-1} =
$$
\be
\label{nr32}
=  \sum\limits_{t=1}^N \sum\limits_{s=1}^{[{\frac t2}]}\;  s \; \tilde{P}_N(s,t)
\ee
We could deduce the above formula (\ref{nr32}) by noting that the index related to the r.v.
$\numb$ is $s$.
Analogously we find the formula for the mean of the r.v. $\leng$
$$E_N(\leng) = E_N(\numb)+E_N(n_{br})$$
and
$$E_N(n_{br})={\partial \over {\partial w}} \zn\uvw \uno=$$
$$
\left ( {2 \over 3}\right )^{N-1}  \sum\limits_{t=1}^N \sum\limits_{s=1}^{[{\frac t2}]}\; {{t-2s} \over {2^{t-s}}}  {N-t+s \choose s} {t-s-1 \choose s-1}= $$
$$
  \sum\limits_{t=1}^N \sum\limits_{s=1}^{[{\frac t2}]}\;(t-2s) \tilde{P}_N(s,t)$$
Hence
\be
\label{nr33}
 E_N(\leng) =
\left ( {2 \over 3}\right )^{N-1}  \sum\limits_{t=1}^N \sum\limits_{s=1}^{[{\frac t2}]}\; {{t-s} \over {2^{t-s}}}  {N-t+s \choose s} {t-s-1 \choose s-1}
\ee
In the next proposition we prove a recursive asymptotic formula, that links the mean of the r.v. $\numb$ to that one of $\leng$, for $N$ large enough.

\begin{proposition}
The following asymptotics
\be
\label{nr34}
E_N(\numb) \asymp {N \over 3} - E_N(\leng)  + {\frac 23} \; E_{N-1}(\leng)
\ee
holds, for $N$ large enough.
\end{proposition}

\begin{remark}
From Corollary 2.3 we have that under $P_N$, $\leng \sim B(N-1, \frac 13)$ (binomial distribution), so that under
$P_{N-1}$, $\leng \sim B(N-2, \frac 13)$. Hence
\bea
\label{nr35}
E_N(\leng) = {{N-1}\over 3} \nonumber\\
E_{N-1}(\leng) = {{N-2}\over 3}
\eea
\end{remark}
From Remark 3.1 the next Corollary follows:

\begin{cor} For large $N$,
\be
\label{nr36}
E_N(\numb) \asymp {{2N-1}\over 9}
\ee
\end{cor}
\nt
{\bf Proof:} Applying the formula (\ref{nr34}) and taking into account (\ref{nr35}) we obtain
$$E_N(\numb)  \asymp {N \over 3} -  {{N-1}\over 3}  + {\frac 23} \;  {{N-2}\over 3}=
{{2N-1}\over 9}$$
\hfill{$\Box$}
\begin{remark}
By identity (\ref{nr1}) and from (\ref{nr35}) and (\ref{nr36}) we are able to calculate asymptotically the single point number's mean. In fact
$$ E_N(\gamma_b+\gamma_r) =N - E_N(\numb) - E_N(\leng) \asymp$$
$$
N - {{2N-1}\over 9} - {{N-1}\over 3} ={\frac 49}(N+1)$$
Note that if we consider only the first order of the asymptotics with respect to $N$, we have
$$E_N(\gamma_b+\gamma_r) \asymp {\frac 49}N \asymp 2 E_N(\numb)$$
that is, for the present model the expected number of single points
is asymptotically twice the expected number of dimers. Moreover,
fixing the number of single points, the conditional probability
distribution of $\gamma_b$ ($\gamma_r$) is binomial and symmetric.
\end{remark}

\medskip
\nt
{\bf Proof of Proposition 3.1}: From (\ref{nr32}) we have
\be
\label{nr37}
E_N(\numb) =
\left ( {2 \over 3}\right )^{N-1}  \sum\limits_{t=1}^N \sum\limits_{s=1}^{[{\frac t2}]}{{N-t+s} \over {2^{t-s}}} {N-t+s-1 \choose s-1} {t-s-1 \choose s-1}
\ee
We applied in (\ref{nr37}) the identity
$$ s {N-t+s \choose s} =
(N-t+s) {N-t+s-1 \choose s-1} $$
Moreover, by Pascal's identity
$${N-t+s-1 \choose s-1}= {N-t+s \choose s} - {N-t+s-1 \choose s}$$
 the mean $E_N(\numb)$ becomes
\be
\label{nr38}
\left ( {2 \over 3}\right )^{N-1}  \sum\limits_{t=1}^N \sum\limits_{s=1}^{[{\frac t2}]} {{N-t+s} \over {2^{t-s}}}  \left [ {N-t+s \choose s} - {N-1-t+s \choose s} \right ] {t-s-1 \choose s-1}
\ee
Taking into account the definition of the probability distribution ${\tilde P}_N$ in  (\ref{nr5}), the formula  (\ref{nr38}) becomes
\be
\label{nr39}
E_N(\numb) =
 \sum\limits_{t=1}^N \sum\limits_{s=1}^{[{\frac t2}]} (N-t+s) \tilde{P}_N(s,t) -
{\frac 23} \sum\limits_{t=1}^{N-1} \sum\limits_{s=1}^{[{\frac t2}]} (N-t+s) \tilde{P}_{N-1}(s,t)
\ee
In fact the normalizing constant for $\tilde{P}_{N-1}$ is $C_{N-1}=\left ({3 \over 2} \right )^{N-2}$, instead of
$C_N$.
The first sum in (\ref{nr39}) gives
$$N \sum\limits_{t=1}^N \sum\limits_{s=1}^{[{\frac t2}]} \tilde{P}_N(s,t) -  \sum\limits_{t=1}^N \sum\limits_{s=1}^{[{\frac t2}]} (t-s) \tilde{P}_N(s,t) =$$
$$N \left [ 1 - \left ( {2 \over 3}\right )^{N-1} \right ] - E_N(\leng) \asymp $$
\be
\label{nr310}
N - E_N(\leng)
\ee
by (\ref{nr5}).

\nt
Analogously for the second sum in (\ref{nr39}) one has
$$- {\frac 23} \left [ N \sum\limits_{t=1}^{N-1} \sum\limits_{s=1}^{[{\frac t2}]} \tilde{P}_{N-1}(s,t) -  \sum\limits_{t=1}^{N-1} \sum\limits_{s=1}^{[{\frac t2}]} (t-s) \tilde{P}_{N-1}(s,t) \right ]=$$
$$- {\frac 23}N \left [ 1 - \left ( {2 \over 3}\right )^{N-2} \right ] + {\frac 23} E_{N-1}(\leng) \asymp$$
\be
\label{nr311}
 - {\frac 23}N + {\frac 23}E_{N-1}(\leng)
\ee
Both formulas (\ref{nr310}) and (\ref{nr311}) give (\ref{nr34}) and Proposition 3.1 is so proved.

\hfill{$\Box$}

In order to find the variance of $\numb$ we need an analogous recursive formula for the second factorial
moment of $\numb$ and that one of $\leng$, whose distribution is well known. By (\ref{nr32}) the second factorial moment of $\leng$ is of the form

$$
E_N[(\numb)(\numb -1)]=$$
\be
\label{nr312}
\sum\limits_{t=1}^N \sum\limits_{s=1}^{[{\frac t2}]} \; s(s-1)\tilde{P}_N(s,t)
\ee
Analogously it easy to see that the second factorial moment of $\leng$ is
 $$E_N[(\leng)(\leng -1)]=$$
 \be
 \label{nr313}
  \sum\limits_{t=1}^N \sum\limits_{s=1}^{[{\frac t2}]} (t-s)(t-s-1)\tilde{P}_N(s,t)
  \ee
  In fact, the indices related to the r.v.s $\numb$ and $\leng$ are $s$ and $t-s$ respectively.

  We generalize the asymptotic recursive formula for the first moments, given in Proposition 3.1 to the second factorial moment of $\numb$ in terms of the first two factorial moments of $\leng$.

\begin{proposition} For $N$ large enough, the following
 asymptotic recursive formula
$$
E_N[(\numb)(\numb-1)] \asymp {{N(N-1)} \over 9} +$$
$$
 - 2(N-1) [E_N(\leng)  - {\frac 43} \; E_{N-1}(\leng) +
{\frac 49} \; E_{N-2}(\leng)] +
$$
$$
E_N[(\leng)(\leng-1)]  - {\frac 43} E_{N-1}[(\leng)(\leng-1)] +
$$
\be
\label{nr314}
 {\frac 49} E_{N-2}[(\leng)(\leng-1)]
\ee
holds.
\end{proposition}
\medskip
\nt
{\bf Proof of Proposition 3.2}: From (\ref{nr312}) and taking into account the proof of Proposition 3.1, we have
$$
E_N[(\numb)(\numb-1)]  =
$$
$$
\left ( {2 \over 3}\right )^{N-1}  \sum\limits_{t=1}^N \sum\limits_{s=1}^{[{\frac t2}]}{{(N-t+s)(N-t+s-1)} \over {2^{t-s}}} {N-t+s-2 \choose s-2} {t-s-1 \choose s-1}
$$
because of the identity
$$ s(s-1) {N-t+s \choose s} =
(N-t+s)(N-t+s-1) {N-t+s-2 \choose s-2} $$
Considering the expressions   (\ref{nr312}) and  (\ref{nr313}) of the second factorial moments of $\numb$ and $\leng$ respectively, we rewrite $(N-t+s)(N-t+s-1)$ in a suitable form
\be
\label{nr315}
(N-t+s)(N-t+s-1)= N(N-1) -2(N-1)(t-s)+(t-s)(t-s-1)
\ee
Moreover, as in Proposition 3.1, we apply (this time twice) Pascal's formula
$${N-t+s-2 \choose s-2}= {N-t+s-1 \choose s-1} - {N-t+s-2 \choose s-1}=$$
\be
\label{nr316}
{N-t+s \choose s} - 2{N-t+s-1 \choose s}+{N-t+s-2 \choose s}
\ee
By (\ref{nr315}) and  (\ref{nr316}) the second factorial moment of $\numb$ becomes
$$
E_N[(\numb)(\numb-1)]  =
$$
$$
\left ( {2 \over 3}\right )^{N-1}  \sum\limits_{t=1}^N \sum\limits_{s=1}^{[{\frac t2}]}
[N(N-1) -2(N-1)(t-s)+(t-s)(t-s-1)] \times
$$
$$
\left [ {N-t+s \choose s} - 2{N-t+s-1 \choose s}+{N-t+s-2 \choose s}\right ] {t-s-1 \choose s-1}
{{1} \over {2^{t-s}}}=$$
$$
 \sum\limits_{t=1}^N \sum\limits_{s=1}^{[{\frac t2}]} [N(N-1) -2(N-1)(t-s)+(t-s)(t-s-1)] \tilde{P}_N(s,t) +
 $$
 $$
- {\frac 43}  \sum\limits_{t=1}^{N-1} \sum\limits_{s=1}^{[{\frac t2}]} [N(N-1) -2(N-1)(t-s)+(t-s)(t-s-1)] \tilde{P}_{N-1}(s,t) +
 $$
\be
\label{nr317}
+{\frac 49}  \sum\limits_{t=1}^{N-2} \sum\limits_{s=1}^{[{\frac t2}]} [N(N-1) -2(N-1)(t-s)+(t-s)(t-s-1)] \tilde{P}_{N-2}(s,t)
\ee

\nt
where the normalizing constant for $\tilde{P}_{N-2}$ is $C_{N-2}=\left ({3 \over 2} \right )^{N-3}$.

\nt
The first sum of the right hand side of (\ref{nr317}) becomes
$$N (N-1) \sum\limits_{t=1}^N \sum\limits_{s=1}^{[{\frac t2}]} \tilde{P}_N(s,t) -  2(N-1) \sum\limits_{t=1}^N \sum\limits_{s=1}^{[{\frac t2}]} (t-s) \tilde{P}_N(s,t) +$$
$$
\sum\limits_{t=1}^N \sum\limits_{s=1}^{[{\frac t2}]} (t-s)(t-s-1) \tilde{P}_N(s,t) =
$$
$$N(N-1) \left [ 1 - \left ( {2 \over 3}\right )^{N-1} \right ] - 2(N-1)E_N(\leng) + $$
$$
E_N[(\leng)(\leng-1)]  \asymp N(N-1)- 2(N-1)E_N(\leng) + $$
\be
\label{nr318}
E_N[(\leng)(\leng-1)]
\ee
\nt
As above the second sum in (\ref{nr317}) is of the form
$$ - {\frac 43}  N (N-1) \left [ 1 - \left ( {2 \over 3}\right )^{N-2} \right ]
+ {\frac 83}(N-1)E_{N-1}(\leng) + $$
$$
- {\frac 43} E_{N-1}[(\leng)(\leng-1)]  \asymp
$$
$$ - {\frac 43} N(N-1) + {\frac 83} (N-1)E_{N-1}(\leng) + $$
\be
\label{nr319}
- {\frac 43} E_{N-1}[(\leng)(\leng-1)]
\ee

\nt
Finally the last sum in (\ref{nr317}) gives the following contribution
$$  {\frac 49}  N (N-1) \left [ 1 - \left ( {2 \over 3}\right )^{N-3} \right ]
- {\frac 89}(N-1)E_{N-2}(\leng) + $$
$$
+{\frac 49} E_{N-2}[(\leng)(\leng-1)]  \asymp
$$
$$ {\frac 49} N(N-1) - {\frac 89} (N-1)E_{N-2}(\leng) + $$
\be
\label{nr320}
 {\frac 49} E_{N-2}[(\leng)(\leng-1)]
\ee
Putting together  (\ref{nr318})-(\ref{nr320}) we get (\ref{nr314}) and Proposition 3.2 is so proved.

\hfill{$\Box$}

Now we are able to calculate the second factorial moment and the variance of $\numb$

\begin{cor} For large $N$,  the following asymptotics
\be
\label{nr321}
E_N[(\numb)(\numb-1)] \asymp {{4}\over 81}(N-1)(N-3)
\ee
holds.
\end{cor}

\nt
{\bf Proof:} Since $\leng \sim B(N-1, \frac 13)$ for sequences $\xi_N$ of length $N$, one can easily find its second factorial moment
\be
\label{nr322}
E_N[(\leng)(\leng-1)] ={{(N-1)(N-2)} \over 9}
\ee
Analogously we have
\bea
E_{N-1}[(\leng)(\leng-1)] ={{(N-2)(N-3)} \over 9}\nonumber \\
\label{nr323}
E_{N-2}[(\leng)(\leng-1)] ={{(N-3)(N-4)} \over 9}\
\eea
From (\ref{nr35}), (\ref{nr322}), (\ref{nr323}) and Proposition 3.2 we obtain the second factorial moment of $\numb$

$$E_N[(\numb)(\numb-1)]  \asymp $$
$$ {{N(N-1)} \over 9}- 2 (N-1) \left[   {{N-1}\over 3} -  {\frac 43}
 {{N-2}\over 3} + {\frac 49}
 {{N-3}\over 3}  \right ]+$$
 $${{(N-1)(N-2)} \over 9}-  {\frac 43} {{(N-2)(N-3)} \over 9} + {\frac 49} {{(N-3)(N-4)} \over 9} =$$
 $$ {{N(N-1)} \over 9}  -  {\frac 23} (N-1)^2 + (N-1)(N-2) - {\frac 8{27}} (N-1)(N-3)+ $$
 $$-
 {\frac 4{27}} (N-2)(N-3) + {\frac 4{81}} (N-3)(N-4)=$$
 $${{N-1} \over 9}[N+9(N-2)-6(N-1)]- {\frac {32}{81}} (N-1)(N-3)=$$
 $$ {{4}\over 81}(N-1)(N-3)$$

 \hfill{$\Box$}

 \begin{cor} For $N$ large enough
\be
\label{nr324}
Var_N(\numb)\asymp {{2}\over 81}(3N+1)
\ee
The symbol $Var_N$ indicates the variance associated to probability measure $P_N$.
\end{cor}

\nt
{\bf Proof:} One has that
$$Var_N(\numb)= E_N[(\numb)(\numb-1)] + E_N[\numb] -(E_N[\numb] )^2\asymp$$
$$  {{4}\over 81}(N-1)(N-3) +  {{2N-1}\over 9} -  {{(2N-1)^2}\over 81} =$$
$${{1}\over 81}[4(N-1)(N-3)+9(2N-1) -(2N-1)^2]=$$
$${{2}\over 81}(3N+1)$$

  \hfill{$\Box$}

\begin{remark}
Note that in the variance formula (\ref{nr324}) for $\numb$ the second order term with respect to $N$ disappears, so that the accuracy of the first order terms with respect to $N$ is important. Nevertheless, from the asymptotics for the mean and the variance of $\numb$ ((\ref{nr36}) and (\ref{nr324})) one can deduce that the distribution of $\numb$ is asymptotically not binomial. In fact $E_N(\numb) \asymp {2 \over 9} N$ and $Var_N(\numb) \asymp  {2 \over 27} N \neq  {2 \over 9}
\cdot  {7 \over 9} N$. In the next section we prove that it is asymptotically gaussian, for large $N$, i.e. a C.L.T. holds.
\end{remark}

\section{Central Limit Theorem for dimers' number}
 \setcounter{secnum}{\value{section}
 \setcounter{equation}{0}
 \renewcommand{\theequation}{\mbox{\arabic{secnum}.\arabic{equation}}}}

In the present section we study the asymptotic distribution of the dimers' number, in particular we prove a C.L.T. for the joint probability distribution of the total number of dimers and single points ($\tnumb$), analyzed in Section 2, and the number of dimers ($\numb$). The limit distribution is a bivariate gaussian distribution with correlation coefficients equal to $- {1 \over {\sqrt{3}}}$. The proof is a generalization of De Moivre-Laplace's Theorem.

\begin {theorem}
A C.L.T. holds for the joint probability distribution
$$
P_N(\numb=s; \tnumb=k)=
$$
\be
\label{nr41}
 {1 \over {\sqrt{(2 \pi)^2 \left ( {\frac 2{27}} N \right )  \left ( {\frac 2{9}} N \right )  \left ( {\frac 2{3}}  \right ) }}}
e^{-{\frac 34} \left (  x^2 + {{2\sqrt{3}} \over 3 } xy + y^2 \right )}
 \left (1 + r_{x,y}(N)) \right )
 \ee
where
\be
\label{nr42}
 x={{ s - {\frac 29}N} \over { {{\sqrt{6N}}\over 9 }  }} \qquad
y={{ k - {\frac 23}N} \over {{{\sqrt{2N}}\over 3 } }}
\ee
Moreover $\lim_{N \to \infty} r_{x,y}(N) =0$, uniformly with respect to $x$ and $y$, defined in (\ref{nr42}),
belonging to finite intervals $(-C, C)$ and $(-A, A)$ respectively, with $C$ and $A$ positive real  constants.
\end{theorem}

\nt
{\bf Proof:} As in Section 2 we perform the change of variables $k=N-t+s$ and $s=s$ on ${\tilde{P}}_N(s,t)$, defined in (\ref{nr5}). The index $k$ indicates the total number of dimers and single points. Then the probability ${\tilde{P}}_N(s,t)$ becomes
$$
P_N(\numb=s; \tnumb=k)= {\tilde{P}}_N(s,N-k+s)=$$
\be
\label{nr43}
{k \choose s } {N-k-1 \choose s-1} \left ({2 \over 3} \right )^{k-1}
 \left ({1 \over 3} \right )^{N-k}
\ee
Taking into account that the indices $s$ and $k$ are both of order $N$, we can forget
 in (\ref{nr43}) the constants, i.e. $k-1 \asymp k$ and $N-k-1 \asymp N-k$, as $N \to \infty$.

 As in De Moivre-Laplace's Theorem, we apply Stirling's formula to the binomial coefficients in   (\ref{nr43}). In the present model we have two binomial coefficients instead of one, so that the calculus becomes heavier than in De Moivre-Laplace's Theorem. We write
 $$
 {k \choose s } {N-k \choose s}  \left ({2 \over 3} \right )^{k}
 \left ({1 \over 3} \right )^{N-k}= $$
  $$
 {{k!} \over {s! (k-s)!}} {{(N-k)!} \over {s!(N-k-s)!}}  \left ({2 \over 3} \right )^{k}
 \left ({1 \over 3} \right )^{N-k}= $$
 $$
   {{k^k} \over {s^s (k-s)^{k-s}}} {{(N-k)^{N-k}} \over {s^s(N-k-s)^{N-k-s}}}  \left ({2 \over 3} \right )^{k}
 \left ({1 \over 3} \right )^{N-k} \cdot
 $$
 \be
 \label{nr44}
{{e^{\lambda_{N,s,k}}}\over{\sqrt{2 \pi {{s(k-s)}\over{k}}}\sqrt{2 \pi {{s(N-k-s)}\over{N-k}}}}}
 \ee
 where
 $${\lambda}_{N,s,k}\equiv {\lambda}_{k}-{\lambda}_{s}- {\lambda}_{k-s} + {\lambda}_{N-k}- {\lambda}_{s}-{\lambda}_{N-k-s}$$
 and
 $${1 \over {12n+1}} \leq \lambda_{n} \leq  {1 \over {12n}} , \; \textnormal{for any $n \in \mathbb{N}$}$$
 Then we rewrite (\ref{nr44}) as
 $$
 \underbrace{e^{\lambda_{N,s,k}}}_{A_{N,s,k}} \cdot
\underbrace{ {{1}\over{\sqrt{2 \pi {{s(k-s)}\over{k}}}\sqrt{2 \pi {{s(N-k-s)}\over{N-k}}}}}}_{B_{N,s,k}} \cdot
 $$
 \be
 \label{nr45}
\underbrace{ \left ( {k \over{3s}}\right )^s  \left ( {{2k} \over{3(k-s)}}\right )^{k-s}
 \left ( {{2(N-k)} \over{3s}}\right )^{s}  \left ( {{N-k} \over{3(N-k-s)}}\right )^{N-k-s} }_{C_{N,s,k}}
 \ee
 In order to find an asymptotics for (\ref{nr45}) we recall that the r.v. $\tnumb$ has binomial distribution $B(N-1, \frac 23)$ (Corollary 2.1) and hence
 $E_N(\tnumb) \asymp {\frac 23} (N-1) \asymp  {\frac 23} N$ and
 $Var_N(\tnumb) \asymp {\frac 29} (N-1) \asymp  {\frac 29} N$. Moreover in the previous section we proved that $E_N(\numb) \asymp {{2N-1}\over 9}  \asymp  {\frac 29} N$ and that
 $Var_N(\numb) \asymp {\frac 2{81}} (3N+1) \asymp  {\frac 6{81}} N$.

 According to De Moivre-Laplace's Theorem we normalize the variables $\tnumb$ and $\numb$, by using the asymptotics of the moments and performing then the change of variables
 (\ref{nr42}) so that
 \be
 \label{nr46}
  \left \{ \begin{array} {ll}
{s={\frac 29}N+ {{\sqrt{6N}}\over 9}x}  & \textnormal{for  $x \in (-C,C)$} \\
{k={\frac 23}N+ {{\sqrt{2N}}\over 3}y}
 & \textnormal{for  $y \in (-A,A)$}
  \end{array} \right.
 \ee
 Taking into account (\ref{nr46}), we consider the first factor $A_{N,s,k}$ in (\ref{nr45}) with respect to the variables $x$ and $y$, in particular
 $$\lambda_{N,s,k} \leq
 $$
 $${1 \over {12\left ( {\frac 23} N + {{\sqrt{2N}}\over 3}y \right )}} -
 {2 \over {12\left ( {\frac 29} N + {{\sqrt{6N}}\over 9}x  \right )+1 }}  - {1 \over {12\left ( {\frac 49} N - {{\sqrt{6N}}\over 9}x  + {{\sqrt{2N}}\over 3}y \right )+1 }}$$
 \be
 \label{nr47}
 +   {1 \over {12\left ( {\frac N3}  -  {{\sqrt{2N}}\over 3}y \right ) }} -
   {1 \over {12\left ( {\frac N9}  - {{\sqrt{6N}}\over 9}x  - {{\sqrt{2N}}\over 3}y \right )+1 }}
 \ee
 Since $x,y$ belong to bounded intervals one can estimate $\lambda_{N,s,k}$ uniformly from above with respect to $x$ and $y$ and we get a bound of order ${1 \over N}$ that doesn't depend on $x$ and $y$. Analogously, we can find a uniform lower bound of order ${1 \over N}$.
Therefore
 $$A_{N,s,k}=1+ r^{(1)}_{x,y}(N)$$
with $r^{(1)}_{x,y}(N)\rightarrow 0$, as $N \to \infty$, uniformly with respect to $x$ and $y$.

\nt
The second factor in  (\ref{nr45}) is $B_{N,s,k}$
 $$B_{N,s,k}=  {{1}\over{\sqrt{2 \pi {{s(k-s)}\over{k}}}\sqrt{2 \pi {{s(N-k-s)}\over{N-k}}}}}$$
 From (\ref{nr46}) we estimate
 $${{s(k-s)}\over{k}} = {{
 \left ( {\frac 29} N + {{\sqrt{6N}}\over 9}x  \right )
 \left ( {\frac 49} N - {{\sqrt{6N}}\over 9}x  + {{\sqrt{2N}}\over 3}y \right )
 }\over{\left ( {\frac 23} N + {{\sqrt{2N}}\over 3}y \right )}
 }=$$
 $$ {4 \over {27}}\; N
 { {{\left ( 1 + {{\sqrt{6 \over N}} }\;{x \over 2}  \right )\left ( 1 - {{\sqrt{6 \over N}} }\;{x \over 4}
  +{\frac 34}{{\sqrt{2 \over N}} }\;{y}\right )}}\over {\left ( 1 + {{\sqrt{2 \over N}} }\;{y \over 2}  \right )}}=$$
  \be
  \label{nr48}
  {4 \over {27}}\; N \left ( 1 + r^{(2)}_{x,y}(N)\right )
  \ee
with $r^{(2)}_{x,y}(N)\rightarrow 0$, as $N \to \infty$, uniformly with respect to $x$ and $y$, since $x \in (-C, C)$ and
 $y \in (-A, A)$, so that the upper and lower bounds of $B_{N,s,k}$
 do not depend on $x$ and $y$.

\nt
Analogously
 $${{s(N-k-s)}\over{N-k}} =  {2 \over {27}}N \left ( 1 + r^{(3)}_{x,y}(N)\right )$$
with $r^{(3)}_{x,y}(N)\rightarrow 0$, as $N \to \infty$, uniformly  with respect to $x$ and $y$.
  Therefore
  \be
  \label{nr49}
   B_{N,s,k}\asymp {1\over {2 \pi {{\sqrt{6N}}\over 9} {{\sqrt{2N}}\over 3}
  \sqrt{ {{{2}}\over 3}}
   }}
    \ee
   as $N \to \infty$.

\nt
Note that
   $$
    B_{N,s,k}\asymp {1\over {2 \pi \sqrt{Var_N(\numb)}  \sqrt{Var_N(\tnumb)} \sqrt{1 - \rho^2}}}
  $$
where $\rho = \pm {1 \over {\sqrt{3}}}$ is the correlation coefficient, whose sign will be determined later.

\nt
Finally we consider the logarithm of the last factor $C_{N,s,k}$ in  (\ref{nr45})
$$
\ln{C_{N,s,k}} =
- s \ln{ \left ( {{3s} \over{k}}\right )}  -(k-s)  \ln{ \left ( {{3(k-s)} \over{2k}}\right )}  +$$
$$
- s \ln{ \left ( {{3s} \over{2(N-k)}}\right )}  -(N-k-s)  \ln{ \left ( {{3(N-k-s)} \over{N-k}}\right )}  \equiv$$
\be
\label{nr410}
\sum_{i=1}^4{C^i_{N,s,k}}
\ee
We express now each term of the sum in  (\ref{nr410}) ${C^i_{N,s,k}}$, $i=1,2,3,4$ in terms of $x$ and $y$, defined in  (\ref{nr42}). We start with ${C^1_{N,s,k}}$
$${C^1_{N,s,k}} \equiv
- s \ln{ \left ( {{3s} \over{k}}\right )} = $$
$$
- \left ( {\frac 29} N + {{\sqrt{6N}}\over 9}x  \right ) \ln{\left (
{
{{\frac 23} N + {{\sqrt{6N}}\over 3}x}
\over
{{\frac 23} N + {{\sqrt{2N}}\over 3}y}
}
\right )}=
$$
$$
- {{\sqrt{2N}} \over 9} \left ( \sqrt{2N} + \sqrt{3} \; x   \right ) \ln{\left (1 +  {{\sqrt{3} \; x - y}\over {\sqrt{2N} + y}} \right )}
$$
Since the last logarithm above is of the form $\ln(1+z)$, with $z \to 0$, we can expand it around $z=0$, $\ln(1+z) = z - {{z^2} \over 2} + o(z^2)$, as $z \to 0$. The same is true for each logarithm function present in any ${C^i_{N,s,k}}$, $i=1,2,3,4$. So
${C^1_{N,s,k}} $ becomes, as $N \to \infty$,
\be
\label{nr411}
{C^1_{N,s,k}} \asymp - {{\sqrt{2N}} \over {18}} \;
 {{\left (\sqrt{2N} + \sqrt{3} \; x)(\sqrt{3} \; x - y \right )\left (2 \sqrt{2N} - \sqrt{3} \; x +3 y \right )}\over {\left (\sqrt{2N} + y \right )^2}}
\ee

\nt
Analogously for ${C^2_{N,s,k}}$
$${C^2_{N,s,k}} \equiv
- (k-s) \ln{ \left ( {{3(k-s)} \over{2k}}\right )} = $$
$$
- \left ( {\frac 49} N - {{\sqrt{6N}}\over 9}x + {{\sqrt{2N}}\over 3}y   \right )
 \ln{\left (
{
{{\frac 43} N - {{\sqrt{6N}}\over 3}x +{{\sqrt{2N}}}\; y}
\over
{{\frac 43} N +2{{\sqrt{2N}}\over 3}y}
}
\right )}=
$$
$$
- {{\sqrt{2N}} \over 9} \left ( 2 \sqrt{2N} - \sqrt{3} \; x +3 y   \right ) \ln{\left (1 +  {{- \sqrt{3} \; x + y}\over {2(\sqrt{2N} + y)}} \right )} \asymp
$$
\be
\label{nr412}
+ {{\sqrt{2N}} \over {72}}  \;
 {{\left ( 2 \sqrt{2N} - \sqrt{3} \; x +3 y   \right )\left (\sqrt{3} \; x - y \right ) \left (4 \sqrt{2N} + \sqrt{3} \; x +3 y \right )}\over {\left (\sqrt{2N} + y \right )^2}}\ee

Then ${C^3_{N,s,k}}$
$${C^3_{N,s,k}} \equiv
-s \ln{ \left ( {{3s} \over{2(N-k)}}\right )} = $$
$$
- \left ( {\frac 29} N + {{\sqrt{6N}}\over 9}x    \right )
 \ln{\left (
{
{2N + {{\sqrt{6N}}}x}
\over
{2N - 2{{\sqrt{2N}}}y}
}
\right )}=
$$
$$
- {{\sqrt{2N}} \over 9} \left (  \sqrt{2N} + \sqrt{3} \; x    \right ) \ln{\left (1 +  {{ \sqrt{3} \; x + 2 y}\over {\sqrt{2N} -  2 y}} \right )} \asymp
$$
\be
\label{nr413}
- {{\sqrt{2N}} \over {18}} \;
 {{\left (  \sqrt{2N} + \sqrt{3} \; x    \right )\left (\sqrt{3} \; x + 2 y \right ) \left (2 \sqrt{2N} - \sqrt{3} \; x -6 y \right )}\over {\left (\sqrt{2N} - 2 y \right )^2}}\ee

 It remains to see ${C^4_{N,s,k}}$
$${C^4_{N,s,k}} \equiv
-(N-k-s) \ln{ \left ( {{3(N-k-s)} \over{N-k}}\right )} = $$
$$
- \left ( {\frac N9}  - {{\sqrt{6N}}\over 9}x   -  {{\sqrt{2N}}\over 3}y  \right )
 \ln{\left (
{
{{\frac N3} - {{\sqrt{6N}}\over 3}x-  {{\sqrt{2N}}}\; y }
\over
{{\frac N3} - {{\sqrt{2N}}\over 3}y}
}
\right )}=
$$
$$
- {{\sqrt{2N}} \over 9} \left (  \sqrt{N \over 2} - \sqrt{3} \; x  -3 y  \right ) \ln{\left (1 -  {{ \sqrt{3} \; x + 2 y}\over {\sqrt{N \over 2} -   y}} \right )} \asymp
$$
\be
\label{nr414}
 2{{\sqrt{2N}} \over {9}} \;
 {{\left (  \sqrt{N \over 2} - \sqrt{3} \; x -3y    \right )\left (\sqrt{3} \; x + 2 y \right ) \left ( \sqrt{2N} + \sqrt{3} \; x  \right )}\over {\left (\sqrt{2N} - 2 y \right )^2}}\ee

\nt
Summing the last term in (\ref{nr411}) (${C^1_{N,s,k}}$) with that one in
 (\ref{nr412}) (${C^2_{N,s,k}}$) we obtain
 \be
 \label{nr415}
 {C^1_{N,s,k}}+{C^2_{N,s,k}} \asymp
 - {{\sqrt{2N}} \over {24}} \;
 {{\left (2 \sqrt{2N} - \sqrt{3} \; x+3y)(\sqrt{3} \; x - y \right )^2}\over {\left (\sqrt{2N} + y \right )^2}}
 \ee

\nt
Summing the last term in (\ref{nr413}) (${C^3_{N,s,k}}$) with that one in
 (\ref{nr414}) (${C^4_{N,s,k}}$) we obtain
 \be
 \label{nr416}
 {C^3_{N,s,k}}+{C^4_{N,s,k}} \asymp
 - {{\sqrt{2N}} \over {6}} \;
 {{\left (\sqrt{2N} + \sqrt{3} \; x)(\sqrt{3} \; x +2 y \right )^2}\over {\left (\sqrt{2N} -2 y \right )^2}}
 \ee

\nt
Finally, the main contribution of (\ref{nr415}) and (\ref{nr416}) is
 $$
  \sum_{i=1}^4{C^i_{N,s,k}} \asymp
  - {{1} \over {12}} \;
 \left [\left (\sqrt{3} \; x - y \right )^2 + 2 \left  (\sqrt{3} \; x + 2 y \right )^2  \right ] =
 $$
 \be
 \label{nr417}
 - {{3} \over {4}} \;
 \left [x^2 + 2{{\sqrt{3}} \over 3} xy +y^2  \right ]
 \ee
\begin{remark}
Note that the last term in (\ref{nr417}) if of the form
$$
-{1 \over {{2(1- \rho^2)}}} \left (x^2 - 2 \rho xy +y^2 \right )
$$
with $\rho= - {1 \over {\sqrt{3}}}$.  We get so a bivariate gaussian distribution with correlation coefficient equal to $ - {1 \over {\sqrt{3}}}$, i.e., the r.v. $\numb$ and $\tnumb$ are negatively correlated.
\end{remark}

\nt
From the previous theorem we can deduce the following result.

\begin{cor} For large $N$, a C.L.T. for the r.v. $\numb$ holds, i.e.
\be
\label{nr418}
P(\numb=s) \asymp {1 \over {\sqrt{ 2 \pi {\frac 2{27}}N}}}
e^{-{{\left (  s - {\frac 29}N \right )^2} \over {{\frac 4{27}}N}}
}
 \ee
for any $s$ defined in (\ref{nr46}).
\end{cor}

\nt
{\bf Proof}: From (\ref{nr43}) we have
$$
P_N(\numb=s) = \sum_{k=s}^{N-s} {k \choose s} { N-k-1 \choose s-1}
\left ({2 \over 3} \right )^{k-1} \left ({1 \over 3} \right )^{N-k}=
$$
\be
\label{nr419}
\left ( \sum\limits_{k=\frac 29N+ {{\sqrt{6N}}\over 9}x}^{\frac 23N- {{\sqrt{2N}}\over 3}A}+
\sum\limits_{k=\frac 23N- {{\sqrt{2N}}\over 3}A}^{\frac 23N+ {{\sqrt{2N}}\over 3}A}+
\sum\limits_{k=\frac 23N+ {{\sqrt{2N}}\over 3}A}^{\frac 79N- {{\sqrt{6N}}\over 9}x} \right )
{\tilde P}_N(s, N-k+s) \\[1ex]
\ee
In (\ref{nr419}), let us denote by $\Sigma_1(A,N), \Sigma_2(A,N)$ and $\Sigma_3(A,N)$ the first, second and third sum, respectively. Using Theorem 4.1
and the integral C.L.T. of De Moivre-Laplace \cite{Gne} one finds for every $\epsilon > 0$ an $N(\epsilon)$ such that
\be
\label{nr420a}
\left|\Sigma_i(A,N)- I_i(A) \right|<\frac{\epsilon}{3},\quad i = 1,2,3
\ee
for any $N > N(\epsilon)$.
The terms $I_i(A)$ are given by
$$
{1 \over {\sqrt{ 2 \pi {\frac 2{27}}N}}}
e^{-{{x^2} \over {2}}}{1 \over {\sqrt{2 \pi}}} \int\limits_{i} e^{-{{y^2} \over{2}}} dy
$$
and the boundaries of the integrals $\int_i$ are fixed as $\int_1=\int_{-\infty}^{-A}$, $\int_2=\int_{-A}^A$ and $\int_3=\int_A^{\infty}$, respectively.
Note that (\ref{nr420a}) proves the statement. \hfill{$\Box$}

\end{document}